\documentclass[a4paper,11pt]{article}
\usepackage{amssymb,amsmath,latexsym,epsf,graphicx}

\newtheorem{thm}{Theorem}[section]

\newtheorem{rem}[thm]{Remark}

\title{The Ring of Support-Classes of $\mathrm{SL}_2(\mathbb F_q)$}
\author{Roland Bacher\footnote{This work has been partially supported by the LabEx PERSYVAL-Lab (ANR--11-LABX-0025). The author is a member of the project-team GALOIS supported by this LabEx.}}

\begin{document}
\maketitle

{\sl Abstract\footnote{
Math. class: Primary: 20G40 , Secondary: 05E30}: We introduce and study a subring $\mathcal{SC}$ of $\mathbb Z[\mathrm{SL}_2(\mathbb F_q)]$ obtained by summing elements of $\mathrm{SL}_2(\mathbb F_q)$
according to their support. The ring $\mathcal SC$ can be used 
for the construction of several association schemes.}



\section{Main results}

Summing elements of the finite group $\mathrm{SL}_2(\mathbb F_q)$
according to their support (locations of non-zero matrix coefficients),
we get seven 
elements (six when working over $\mathbb F_2$) in the integral group-ring 
$\mathbb Z[\mathrm{SL}_2(\mathbb F_q)]$.

Integral linear combinations of these seven elements
form a subring $\mathcal{SC}$, called the \emph{ring of support classes}, of the integral group-ring 
$\mathbb Z[\mathrm{SL}_2(\mathbb F_q)]$. Supposing $q>2$, we get 
thus a $7-$dimensional algebra 
$\mathcal{SC}_{\mathbb K}=
\mathcal{SC}\otimes_{\mathbb Z}\mathbb K$ over a field $\mathbb K$
when considering $\mathbb K-$linear combinations.

This paper is devoted to the definition
and the study of a few features of $\mathcal{SC}$.

More precisely, in Section \ref{sectringsuppcl} we prove that 
the ring of support-classes $\mathcal{SC}$
is indeed a ring by computing its structure-constants.

Section \ref{sectalgprop} describes the structure
of $\mathcal{SC}_{\mathbb Q}=\mathcal{SC}\otimes_{\mathbb Z}\mathbb Q$ as a 
semi-simple algebra independent of $q$ for $q>2$.

In Section \ref{sectass} we recall the definition of association schemes 
and use $\mathcal{SC}$ for the construction of
hopefully interesting examples.

Finally, we study a few representation-theoretic aspects in 
Section \ref{sectrepr}.


\section{The ring of support-classes}\label{sectringsuppcl}

Given subsets $\mathcal A,\mathcal B,\mathcal C,\mathcal D$ of a finite field $\mathbb F_q$,
we denote by
$$\left(\begin{array}{cc}\mathcal A&\mathcal B\\\mathcal C&\mathcal D\end{array}\right)=\sum_{(a,b,c,d)\in 
\mathcal A\times \mathcal B\times \mathcal C\times \mathcal D,ad-bc=1}\left(\begin{array}{cc}a&b\\c&d
\end{array}\right)$$
the element of $\mathbb Z[\mathrm{SL}_2(\mathbb F_q)]$
obtained by summing all matrices of $\mathrm{SL}_2(\mathbb F_q)$
with coefficients $a\in \mathcal A,b\in \mathcal B,c\in \mathcal C$ and $d\in \mathcal D$.

Identifying $0$ with the singleton subset $\{0\}$ of $\mathbb F_q$ and 
denoting by $\mathbb F_q^*=\mathbb F_q\setminus\{0\}$ 
the set of all units in $\mathbb F_q$, we consider the seven elements
\begin{eqnarray*}
&&A=\left(\begin{array}{cc}\mathbb F_q^*&0\\0&\mathbb F_q^*\end{array}\right),\
B=\left(\begin{array}{cc}0&\mathbb F_q^*\\\mathbb F_q^*&0\end{array}\right),\ 
C=\left(\begin{array}{cc}\mathbb F_q^*&\mathbb F_q^*\\\mathbb F_q^*&\mathbb F_q^*\end{array}\right),\\
&&D_+=\left(\begin{array}{cc}\mathbb F_q^*&\mathbb F_q^*\\0&\mathbb F_q^*\end{array}\right),\ D_-=
\left(\begin{array}{cc}\mathbb F_q^*&0\\\mathbb F_q^*&\mathbb F_q^*\end{array}\right),\\
&&E_+=
\left(\begin{array}{cc}\mathbb F_q^*&\mathbb F_q^*\\\mathbb F_q^*&0\end{array}\right),\ E_-=\left(\begin{array}{cc}0&\mathbb F_q^*\\\mathbb F_q^*&\mathbb F_q^*\end{array}\right).
\end{eqnarray*}
corresponding to all possible supports of matrices in 
$\mathrm{SL}_2[\mathbb F_q]$.
The element $C$ is of course missing (and the remaining elements 
consist simply of all six matrices in $\mathrm{SL}_2(\mathbb F_2)$)
over $\mathbb F_2$.
For the sake of concision, 
we will always assume that $q$ has more than $2$ elements in the sequel
(there is however nothing wrong with finite fields of characteristic $2$
having at least $4$ elements).

We denote by 
$$\mathcal{SC}=\mathbb Z A+\mathbb Z B+\mathbb ZC+\mathbb ZD_++\mathbb ZD_-
+\mathbb Z E_++\mathbb Z E_-$$
the free $\mathbb Z-$module of rank seven spanned by these seven elements.

We call $\mathcal{SC}$ the \emph{ring of support-classes of 
$\mathrm{SL}_2(\mathbb F_q)$}, a terminology motivated by 
our main result:

\begin{thm}\label{mainthma} $\mathcal{SC}$ is a  
subring of the integral group-ring $\mathbb Z[\mathrm{SL}_2(\mathbb F_q)]$.
\end{thm}

The construction of $\mathcal{SC}$ can be carried over to 
the projective special groups $\mathrm{PSL}_2(\mathbb F_q)$
without difficulties by dividing all structure-constants by $2$
if $q$ is odd. The obvious modifications are left to the reader.

Another obvious variation is to work with matrices in 
$\mathrm{GL}_2(\mathbb F_q)$. This multiplies all 
structure-constants by $(q-1)$ (respectively by $m$
if working with the subgroup of matrices in $\mathrm{GL}_2(\mathbb F_q)$
having their determinants in a fixed multiplicative 
subgroup $M\subset \mathbb F_q^*$  with $m$ elements).

We hope to address a few other variations of our main construction
in a future paper.

Products among generators of $\mathcal{SC}$ are given by 
\begin{eqnarray*}
AX&=&XA=(q-1)X\hbox{ for }X\in\{A,B,C,D_\pm,E_\pm\},\\
B^2&=&(q-1)A,\\ 
BC=CB&=&(q-1)C,\\
BD_+=D_-B&=&(q-1)E_-,\\ 
BD_-=D_+B&=&(q-1)E_+,\\
BE_+=E_-B&=&(q-1)D_-,\\ 
BE_-=E_+B&=&(q-1)D_+,\\
C^2&=&(q-1)^2(q-2)(A+B)+(q-1)(q-3)(q-4)C\\
&&\ +(q-1)(q-2)(q-3)(D_++D_-+E_++E_-),\\
CD_+=CE_-&=&(q-1)(q-3)C+(q-1)(q-2)(D_-+E_+),\\
CD_-=CE_+&=&(q-1)(q-3)C+(q-1)(q-2)(D_++E_-),\\
D_+C=E_+C&=&(q-1)(q-3)C+(q-1)(q-2)(D_-+E_-),\\
D_-C=E_-C&=&(q-1)(q-3)C+(q-1)(q-2)(D_++E_+),\\
D_+^2=E_+E_-&=&(q-1)^2A+(q-1)(q-2)D_+,\\
D_+D_-=E_+^2&=&(q-1)(C+E_-),\\
D_-D_+=E_-^2&=&(q-1)(C+E_+),\\
D_+E_+=E_+D_-&=&(q-1)^2B+(q-1)(q-2)E_+,\\
E_+D_+=D_+E_-&=&(q-1)(C+D_-),\\
E_-D_+=D_-E_-&=&(q-1)^2B+(q-1)(q-2)E_-,\\
D_-^2=E_-E_+&=&(q-1)^2A+(q-1)(q-2)D_-,\\
D_-E_+=E_-D_-&=&(q-1)(C+D_+).\\
\end{eqnarray*}

Easy consistency checks of these formulae are given by the antiautomorphisms
$\sigma$ and $\tau$ obtained respectively by matrix-inversion and matrix-transposition. Their composition $\sigma\circ \tau=\tau\circ \sigma$ is 
of course an involutive automorphism of $\mathbb Z[\mathbb F_q]$ 
which restricts to
an automorphism of $\mathcal{SC}$. It coincides on $\mathcal{SC}$
with the action of the inner automorphism $X\longmapsto 
\left(\begin{array}{cc}0&1\\-1&0\end{array}\right)X
\left(\begin{array}{cc}0&-1\\1&0\end{array}\right)$ of 
$\mathbb Z[\mathrm{SL}_2(\mathbb F_q)]$, fixes $A,B,C$ and 
transposes the elements of the two pairs $\{D_+,D_-\}$ and $\{E_+,E_-\}$.

\begin{rem} The construction of the ring $\mathcal{SC}$ described by Theorem
\ref{mainthma} does not generalise to the matrix-algebra of all 
$2\times 2$ matrices over $\mathbb F_q$. Indeed, 
$\left(\begin{array}{cc}\mathbb F_q^*&0\\\mathbb F_q^*&0\end{array}\right)
\left(\begin{array}{cc}\mathbb F_q^*&\mathbb F_q^*\\0&0\end{array}\right)$
equals, up to a factor $(q-1)$, to the sum of all $(q-1)^3$ possible 
rank $1$ matrices with all four coefficients in $\mathbb F_q^*$.

Square rank one matrices of any size behave however rather well: The set
of all $(2^n-1)^2$ possible sums of rank $1$ matrices of size $n\times n$
with prescribed support is a $\mathbb Z$-basis of a ring (defined by
extending bilinearly the matrix product) after identifying the
zero matrix with $0$.
\end{rem}

\subsection{Proof of Theorem \ref{mainthma}}

We show that the formulae for the products are correct.
We are however not going to prove all $7^2=49$
possible identities but all omitted cases are similar 
and can be derived by symmetry arguments, use of
the antiautomorphisms given by matrix-inversion and 
transposition, or (left/right)-multiplication by
$B$.

Products with $A$ or $B$ are easy and left to the reader.

We start with the easy product $D_+^2$
(the products $$D_+E_+,E_+D_-,E_-D_+,D_-E_-,D_-^2,E_-E_+$$ 
are similar and left to the reader). Since $A+D_+$ is the sum of elements
over the
full group of all $q(q-1)$ unimodular upper-triangular
matrices, we have $(A+D_+)^2=q(q-1)(A+D_+)$ showing that
\begin{eqnarray*}
D_+^2&=&(A+D_+)^2-2AD_+-A^2\\
&=&q(q-1)(A+D_+)-2(q-1)D_+-(q-1)A\\
&=&(q-1)^2A+(q-1)(q-2)D_+.
\end{eqnarray*}

For $D_+D_-$ we consider 
$$\left(\begin{array}{cc}a_1&b_1\\0&1/a_1\end{array}\right)
\left(\begin{array}{cc}a_2&0\\b_2&1/a_2\end{array}\right)
=\left(\begin{array}{cc}a_1a_2+b_1b_2&b_1/a_2\\b_2/a_1&1/(a_1a_2)
\end{array}\right).
$$
Every unimodular matrix $\left(\begin{array}{cc}\alpha&\beta\\
\gamma&\delta\end{array}\right)$
with $\beta,\gamma,\delta$ in $\mathbb F_q^*$ can be realised as 
a summand in the product $D_+D_-$ in exactly $(q-1)$ different 
ways by choosing $a_1$ freely in $\mathbb F_q^*$ and by setting
$$b_1=\frac{\beta}{a_1\delta},a_2=\frac{1}{a_1\delta},b_2=a_1\gamma.$$
This shows $D_+D_-=(q-1)(C+E_-)$. The products
$$E_+^2,D_-D_+,E_-^2,E_+D_+,D_+E_-,D_-E_+,E_-D_-$$
are similar.

In order to compute $D_+C$, we consider $(D_++A)(C+D_-+E_+)$.
Since
$$\left(\begin{array}{cc}
a_1&b_1\\0&1/a_1\end{array}\right)
\left(\begin{array}{cc}
a_2&b_2\\c_2&(1+b_2c_2)/a_2\end{array}\right)
$$
$$\left(\begin{array}{cc}
a_1a_2+b_1c_2&a_1b_2+b_1(1+b_2c_2)/a_2\\
c_2/a_1&(1+b_2c_2)/(a_1a_2)\end{array}\right),$$
every unimodular
matrix $\left(\begin{array}{cc}\alpha&\beta\\\gamma&\delta\end{array}\right)$
with $\gamma\in\mathbb F_q^*$ 
can be realised
in exactly $(q-1)^2$ ways as a summand in $(D_++A)(C+D_-+E_+)$ by 
choosing $a_1,a_2$ freely in $\mathbb F_q^*$ and by setting 
$$b_1=\frac{\alpha-a_1a_2}{a_1\gamma},b_2=\frac{a_1a_2\delta-1}{a_1\gamma},
c_2=a_1\gamma.$$

This shows $(D_++A)(C+D_-+E_+)=(q-1)^2(B+C+D_-+E_++E_-)$. We get thus
\begin{eqnarray*}
&&D_+C\\
&=&(D_++A)(C+D_-+E_+)-A(C+D_-+E_+)-D_+D_--D_+E_+\\
&=&(q-1)^2(B+C+D_-+E_++E_-)-(q-1)(C+D_-+E_+)\\
&&\ -(q-1)(C+E_-)-(q-1)^2B-(q-1)(q-2)E_+\\
&=&(q-1)(q-3)C+(q-1)(q-2)(D_-+E_-).
\end{eqnarray*}
The products
$$CD_+,CE_-,CD_-,CE_+,E_+C,D_-C,E_-C$$
are similar.

Using all previous products, the formula for $C^2$ can now be recovered
from $(A+B+C+D_++D_-+E_++E_-)^2=(q^3-q)(A+B+C+D_++D_-+E_++E_-)$.
We have indeed
\begin{eqnarray*}
C^2&=&(A+B+C+D_++D_-+E_++E_-)^2-\sum_{(X,Y)\not=(C,C)}XY
\end{eqnarray*}
where the sum is over all elements of 
$\{A,B,C,D_+,D_-,E_+,E_-\}^2\setminus(C,C)$. All products of the right-hand-side
are known and determine thus $C^2$.
Equivalently, structure-constants of $C^2$ 
have to be polynomials of degree at most $3$ in $q$.
They can thus also be computed by interpolating the coefficients
in $4$ explicit examples. (Using divisibility by $q-1$, 
computing $3$ examples is in fact enough.)

The existence of these formulae proves Theorem \ref{mainthma}.\hfill$\Box$

\subsection{Matrices for left-multiplication by generators}

Left-multiplications by generators 
with respect to the basis $A,B,C,D_+,D_-,E_+,E_-$
of $\mathcal{SC}$ are encoded by the matrices 
$$M_B=(q-1)\left(\begin{array}{ccccccc}
0&1&0&0&0&0&0\\
1&0&0&0&0&0&0\\
0&0&1&0&0&0&0\\
0&0&0&0&0&0&1\\
0&0&0&0&0&1&0\\
0&0&0&0&1&0&0\\
0&0&0&1&0&0&0\\
\end{array}\right)$$
$$M_C=(q-1)
\left(\begin{array}{ccccccc}
0&0&(q-1)(q-2)&0&0&0&0\\
0&0&(q-1)(q-2)&0&0&0&0\\
1&1&(q-3)(q-4)&q-3&q-3&q-3&q-3\\
0&0&(q-2)(q-3)&0&q-2&q-2&0\\
0&0&(q-2)(q-3)&q-2&0&0&q-2\\
0&0&(q-2)(q-3)&q-2&0&0&q-2\\
0&0&(q-2)(q-3)&0&q-2&q-2&0\\
\end{array}\right)$$

$$M_{D_+}=(q-1)\left(\begin{array}{ccccccc}
0&0&0&q-1&0&0&0\\
0&0&0&0&0&q-1&0\\
0&0&q-3&0&1&0&1\\
1&0&0&q-2&0&0&0\\
0&0&q-2&0&0&0&1\\
0&1&0&0&0&q-2&0\\
0&0&q-2&0&1&0&0\\
\end{array}\right)$$

$$M_{E_+}=(q-1)\left(\begin{array}{ccccccc}
0&0&0&0&0&0&q-1\\
0&0&0&0&q-1&0&0\\
0&0&q-3&1&0&1&0\\
0&1&0&0&0&0&q-2\\
0&0&q-2&1&0&0&0\\
1&0&0&0&q-2&0&0\\
0&0&q-2&0&0&1&0\\
\end{array}\right)$$

The remaining matrices are given by $M_A=\frac{1}{q-1}(M_B)^2$
, $M_{D_-}=\alpha M_{D_+}\alpha$ and $M_{E_-}=
\alpha M_{E_+}\alpha$
where 
$$\alpha=\left(\begin{array}{ccccccc}
1&0&0&0&0&0&0\\
0&1&0&0&0&0&0\\
0&0&1&0&0&0&0\\
0&0&0&0&1&0&0\\
0&0&0&1&0&0&0\\
0&0&0&0&0&0&1\\
0&0&0&0&0&1&0\\
\end{array}\right)$$
is the matrix corresponding to the automorphism $\sigma\circ \tau$.

The map $\{A,B,C,D_+,D_-,E_+,E_-\}\ni X\longmapsto M_X$ extends of course to an isomorphisme between
$\mathcal{SC}$ and 
$$\mathbb ZM_A+\mathbb ZM_B+\mathbb ZM_C+\mathbb ZM_{D_+}+\mathbb ZM_{D_-}
+\mathbb ZM_{E_+}+\mathbb ZM_{E_-}$$
which is a ring. Computations are easier and faster in 
this matrix-ring than in the subring $\mathcal{SC}$
of $\mathbb Z[\mathrm{SL}_2(\mathbb F_q)]$.


\section{Algebraic properties of $\mathcal{SC}_{\mathbb Q}$}\label{sectalgprop}

\subsection{$\mathcal{SC}_{\mathbb Q}$ as a semisimple algebra}\label{subsectsemisimple} 

\begin{thm}\label{thmstructureSC}
The algebra $\mathcal{SC}_{\mathbb Q}=\mathcal{SC}\otimes_{\mathbb Z}\mathbb Q$
is a semi-simple algebra isomorphic to $\mathbb Q\oplus \mathbb Q\oplus \mathbb Q
\oplus M_2(\mathbb Q)$.
\end{thm}

The structure of $\mathcal{SC}_{\mathbb K}=\mathcal{SC}\otimes_{\mathbb Z}\mathbb K$
is of course easy to deduce for any field of characteristic $0$.

The algebra $\mathcal{SC}_{\mathbb K}$ is also semi-simple (and has the same structure) over most finite fields.

Theorem \ref{thmstructureSC} is an easy consequence of the following 
computations:

The center of $\mathcal{SC}_{\mathbb Q}$
has rank $4$. It is spanned by the $3$
minimal central idempotents
\begin{eqnarray*}
\pi_1&=&\frac{1}{q^3-q}(A+B+C+D_++D_-+E_++E_-),\\
\pi_2&=&\frac{q-2}{2(q^2-q)}(A+B)+\frac{1}{q(q-1)^2}C-\frac{q-2}{2q(q-1)^2}
(D_++D_-+E_++E_-),\\
\pi_3&=&\frac{1}{2(q+1)}(A-B)+\frac{1}{2(q^2-1)}(-D_+-
D_-+E_++E_-)
\end{eqnarray*}
and by the central idempotent 
\begin{eqnarray*}
\pi_4&=&\frac{1}{(q+1)(q-1)^2}(2(q-1)A-2C+(q-2)(D_++D_-)-(E_++E_-))
\end{eqnarray*}
which is non-minimal among all idempotents.
The three idempotents $\pi_1,\pi_2,\pi_3$ induce three different 
characters (homomorphisms from $\mathcal{SC}_{\mathbb Q}$ into $\mathbb Q$).
Identifying $1\in \mathbb Q$ with $\pi_i$ in each case,
the three homomorphisms are given by 
$$\begin{array}{||c||c|c|c|c|c||}
\hline\hline
&A&B&C&D_\pm&E_\pm\\
\hline
\pi_1&q-1&q-1&(q-1)^2(q-2)&(q-1)^2&(q-1)^2\\
\pi_2&q-1&q-1&2(q-1)&1-q&1-q\\
\pi_3&q-1&1-q&0&1-q&q-1\\
\hline\hline
\end{array}$$
The idempotent $\pi_1$ is of course simply the augmentation map
counting the number of matrices involved in each generator.

The idempotent $\pi_1+\pi_2+\pi_3+\pi_4=\frac{1}{q-1}A$ is the identity
of $\mathcal{SC}_{\mathbb Q}$.

The idempotent $\pi_4$ projects $\mathcal{SC}_{\mathbb Q}$
homomorphically onto a matrix-algebra of $2\times 2$ matrices.

$\pi_4$ can be written (not uniquely) as a sum of two minimal 
non-central idempotents. We have for example $\pi_4=M_{1,1}+M_{2,2}$
where
\begin{eqnarray*}
M_{1,1}&=&\frac{1}{q^2-1}(A-B)+
\frac{1}{2(q^2-1)}(D_++D_--E_+-E_-),\\
M_{2,2}&=&\frac{1}{q^2-1}(A+B)-\frac{2}{(q+1)(q-1)^2}C+\\
&&\quad +
\frac{q-3}{2(q+1)(q-1)^2}(D_++D_-+E_++E_-).\\
\end{eqnarray*}
Considering also 
\begin{eqnarray*}
M_{1,2}&=&\frac{1}{2(q-1)^2}(D_+-D_-+E_+-E_-),\\
M_{2,1}&=&\frac{1}{2(q^2-1)}(D_+-D_--E_++E_-),\\
\end{eqnarray*}
the elements 
$M_{i,j}$ behave like matrix-units and the map
\begin{eqnarray}\label{exoticisom}
\left(\begin{array}{cc}a&b\\c&d\end{array}\right)\longmapsto
aM_{1,1}+bM_{1,2}+cM_{2,1}+dM_{2,2}
\end{eqnarray}
defines thus an isomorphism from the ring of integral $2\times 2$ matrices
into $\mathcal{SC}_{\mathbb Q}=\mathcal{SC}\otimes_{\mathbb Z}\mathbb Q$.
If $\mathbb F_q$ has odd characteristic, the elements $M_{i,j}$ can be 
realised in $\mathcal{SC}_{\mathbb F_q}=\mathcal{SC}\otimes_{\mathbb Z}\mathbb F_q$.
In particular, Formula (\ref{exoticisom}) gives an \lq\lq exotic\rq\rq, 
non-unital embedding of $\mathrm{SL}_2[\mathbb F_q]$ into
the group-algebra $\mathbb F_q[\mathrm{SL}_2(\mathbb F_q)]$
(in fact, (\ref{exoticisom}) gives an embedding of 
$\mathrm{SL}_2(\mathbb F_r)$ into
$\mathbb F_r[\mathrm{SL}_2(\mathbb F_q)]$ whenever the prime power $r$
is coprime to $2(q^2-1)$).

The nice non-central idempotent $$\pi_3+M_{1,1}=\frac{1}{2(q-1)}(A-B)$$
projects $\mathcal{SC}_{\mathbb Q}$ onto the eigenspace of eigenvalue 
$1-q$ of the map $X\longmapsto BX$.

The projection of $\mathcal{SC}_{\mathbb Q}$ onto the eigenspace of eigenvalue 
$q-1$ of the map $X\longmapsto BX$ is similarly given by
$$\pi_1+\pi_2+M_{2,2}=\frac{1}{2(q-1)}(A+B).$$

\subsection{A few commutative subalgebras of $\mathcal{SC}_{\mathbb Q}$}

The previous section shows that the dimension of a commutative subalgebra of 
$\mathcal{SC}_{\mathbb Q}$ cannot exceed $5$.

The center of $\mathcal{SC}_{\mathbb Q}$ is of course of rank $4$ and
spanned by the four minimal central idempotents $\pi_1,\dots,\pi_4$
of the previous Section.

Splitting $\pi_4=M_{1,1}+M_{2,2}$, we get a maximal commutative subalgebra 
of rank $5$ by considering the vector space spanned by
the minimal idempotents 
$\pi_1,\pi_2,\pi_3,M_{1,1},M_{2,2}$. Equivalently, this vector space
is spanned by $A,B,C,D=D_++D_-,E=E_++E_-$, as shown by the formulae 
for $\pi_i$ and $M_{1,1},M_{2,2}$.

In terms of $A,B,C,D,E$, minimal idempotents $\pi_1,\pi_2,\pi_3,M_{1,1},M_{2,2}$ are given by
\begin{eqnarray*}
\pi_1&=&\frac{1}{q^3-q}(A+B+C+D+E),\\
\pi_2&=&\frac{q-2}{2(q^2-q)}(A+B)+\frac{1}{q(q-1)^2}C-\frac{q-2}{2q(q-1)^2}
(D+E),\\
\pi_3&=&\frac{1}{2(q+1)}(A-B)+\frac{1}{2(q^2-1)}(-D+E),\\
M_{1,1}&=&\frac{1}{q^2-1}(A-B)+
\frac{1}{2(q^2-1)}(D-E),\\
M_{2,2}&=&\frac{1}{q^2-1}(A+B)-\frac{2}{(q+1)(q-1)^2}C+\\
&&\quad +
\frac{q-3}{2(q+1)(q-1)^2}(D+E).\\
\end{eqnarray*}
They define five characters given by 
$\mathbb Q A+\dots+\mathbb Q E\longrightarrow \mathbb Q$ given by
$$\begin{array}{||c||c|c|c|c|c||}
\hline\hline
&A&B&C&D&E\\
\hline
\pi_1&q-1&q-1&(q-1)^2(q-2)&2(q-1)^2&2(q-1)^2\\
\pi_2&q-1&q-1&2(q-1)&2(1-q)&2(1-q)\\
\pi_3&q-1&1-q&0&2(1-q)&2(q-1)\\
M_{1,1}&q-1&1-q&0&(q-1)^2&-(q-1)^2\\
M_{2,2}&q-1&q-1&2(1-q)(q-2)&(q-1)(q-3)&(q-1)(q-3)\\
\hline\hline
\end{array}$$

Moreover,
$A,B,C,F=D+E$ span a commutative $4$-dimensional subalgebra.
Minimal idempotents are 
\begin{eqnarray*}
\pi_1&=&\frac{1}{q^3-q}(A+B+C+F),\\
\pi_2&=&\frac{q-2}{2(q^2-q)}(A+B)+\frac{1}{q(q-1)^2}C-\frac{q-2}{2q(q-1)^2}
F,\\
\pi_3+M_{1,1}&=&\frac{1}{2(q-1)}(A-B),\\
M_{2,2}&=&\frac{1}{q^2-1}(A+B)-\frac{2}{(q+1)(q-1)^2}C+\\
&&\quad +
\frac{q-3}{2(q+1)(q-1)^2}F.\\
\end{eqnarray*}
with character-table
$$\begin{array}{||c||c|c|c|c||}
\hline\hline
&A&B&C&F\\
\hline
\pi_1&q-1&q-1&(q-1)^2(q-2)&4(q-1)^2\\
\pi_2&q-1&q-1&2(q-1)&4(1-q)\\
\pi_3+M_{1,1}&q-1&1-q&0&0\\
M_{2,2}&q-1&q-1&2(1-q)(q-2)&2(q-1)(q-3)\\
\hline\hline
\end{array}$$

$I=A+B,C,F=E+D$ span a commutative $3-$dimensional subalgebra 
of $\mathcal{SC}_{\mathbb C}$.
Generators of this last algebra are sums of elements in 
$\mathrm{SL}_2(\mathbb F_q)$ with supports of given cardinality.
$I$ is the sum of all elements with two non-zero coefficients, 
$C$ contains all elements having only non-zero coefficients and 
$F$ contains all elements with three non-zero coefficients.

Products are given by 
$IX=XI=2(q-1)X$ for $X\in\{I,C,F\}$ 
and 
\begin{eqnarray*}
C^2&=&(q-1)^2(q-2)I+(q-1)(q-3)(q-4)C+(q-1)(q-2)(q-3)F,\\
CF&=&FC=4(q-1)(q-3)C+2(q-1)(q-2)F,\\
F^2&=&4(q-1)^2I+8(q-1)C+2(q-1)^2F.
\end{eqnarray*}

Idempotents are given by 
\begin{eqnarray*}
\pi_1&=&\frac{1}{q^3-q}(I+C+F),\\
\pi_2&=&\frac{q-2}{2(q^2-q)}I+\frac{1}{q(q-1)^2}C-\frac{q-2}{2q(q-1)^2}
F,\\
M_{2,2}&=&\frac{1}{q^2-1}I-\frac{2}{(q+1)(q-1)^2}C+\\
&&\quad +
\frac{q-3}{2(q+1)(q-1)^2}F.\\
\end{eqnarray*}
with character-table
$$\begin{array}{||c||c|c|c||}
\hline\hline
&I&C&F\\
\hline
\pi_1&2(q-1)&(q-1)^2(q-2)&4(q-1)^2\\
\pi_2&2(q-1)&2(q-1)&4(1-q)\\
M_{2,2}&2(q-1)&2(1-q)(q-2)&2(q-1)(q-3)\\
\hline\hline
\end{array}$$

Working over $\mathbb C$ (or over a suitable extension of $\mathbb Q$) and setting 
$$\tilde I=\frac{1}{2(q-1)}I,\tilde C=\frac{1}{\sqrt{2(q-1)^3(q-2)}}C,\tilde F=\frac{1}{\sqrt{8(q-1)^3}}F$$
we get products 
$XY=\sum_{Z\in \{\tilde I,\tilde C,\tilde F\}}N_{X,Y,Z}Z$
defined by symmetric structure-constants $N_{X,Y,Z}=N_{Y,X,Z}=N_{X,Z,Y}$
for all $X,Y,Z\in \{\tilde I,\tilde C,\tilde F\}$.
Up to symmetric permutations, the structure-constants are given by
\begin{eqnarray*}
N_{\tilde I,X,Y}&=&\delta_{X,Y}\\
N_{\tilde C,\tilde C,\tilde C}&=&\frac{(q-3)(q-4)}{\sqrt{2(q-1)(q-2)}},\\
N_{\tilde C,\tilde C,\tilde F}&=&(q-3)\sqrt{\frac{2}{q-1}}\\
N_{\tilde C,\tilde F,\tilde F}&=&\sqrt{\frac{2(q-2)}{q-1}},\\
N_{\tilde F,\tilde F,\tilde F}&=&\sqrt{\frac{q-1}{2}}
\end{eqnarray*}
where $\delta_{X,Y}=1$ if and only if $X=Y$ and $\delta_{X,Y}=0$ otherwise.
The evaluation at $q=3$ leads to particularly
nice structure constants 
with values in $\{0,1\}$.

Algebras with generating systems having symmetric structure-constants
and a character taking
real positive values on generators (satisfied by $\pi_1$ for $q>2$)
are sometimes called
\emph{algebraic fusion-algebras}, see for example \cite{Ba}.

\section{Association schemes and Bose-Mesner algebras}
\label{sectass}

An \emph{association scheme} is a set of $d+1$ square matrices
${\mathcal C}_0,\dots,{\mathcal C}_d$ with coefficients in $\{0,1\}$ such that 
${\mathcal C}_0$ is the identity-matrix, ${\mathcal C}_0+\dots+{\mathcal C}_d$ is the all-one matrix
and $\mathbb Z{\mathcal C}_0+\dots+\mathbb Z{\mathcal C}_d$ is a commutative ring with 
(necessarily integral) structure constants
$p_{i,j}^k=p_{j,i}^k$ defined by ${\mathcal C}_i{\mathcal C}_j={\mathcal C}_j{\mathcal C}_i=\sum_{k=0}^d p_{i,j}^k{\mathcal C}_k$.
An association scheme is \emph{symmetric} if
${\mathcal C}_1,\dots,{\mathcal C}_d$ are symmetric matrices. The algebra (over a field)
generated by the elements ${\mathcal C}_i$ is called a \emph{Bose-Mesner algebra}. See for example the monograph \cite{B} for additional information.

Identifying an element $g$ of $\mathrm{SL}_2(\mathbb F_q)$ with the 
permutation-matrix associated to left-multiplication by $g$
we get a commutative association scheme with $d=5$ if $q\geq 4$ by setting 
$${\mathcal C}_0=\mathrm{Id}, {\mathcal C}_1=A-\mathrm{Id}, {\mathcal C}_2=B,{\mathcal C}_3=C,{\mathcal C}_4=D_++D_-,{\mathcal C}_5=E_++E_-$$
where we consider sums of 
permutation-matrices. All matrices are symmetric. Products with ${\mathcal C}_0,{\mathcal C}_1$
are given by ${\mathcal C}_0X=X{\mathcal C}_0=X, {\mathcal C}_1^2=(q-2){\mathcal C}_0+(q-3){\mathcal C}_1,{\mathcal C}_1Y=Y{\mathcal C}_1=(q-2)Y$
for $X\in\{{\mathcal C}_0,\dots,{\mathcal C}_5\}$ and $Y\in\{{\mathcal C}_2,\dots,{\mathcal C}_5\}$. The remaining products
are given by
\begin{eqnarray*}
{\mathcal C}_2^2&=&(q-1)({\mathcal C}_0+{\mathcal C}_1),\\
{\mathcal C}_2{\mathcal C}_3={\mathcal C}_3{\mathcal C}_2&=&(q-1){\mathcal C}_3\\
{\mathcal C}_2{\mathcal C}_4={\mathcal C}_4{\mathcal C}_2&=&(q-1){\mathcal C}_5\\
{\mathcal C}_2{\mathcal C}_5={\mathcal C}_5{\mathcal C}_2&=&(q-1){\mathcal C}_4\\
{\mathcal C}_3^2&=&(q-1)^2(q-2)({\mathcal C}_0+{\mathcal C}_1+{\mathcal C}_2)+(q-1)(q-3)(q-4){\mathcal C}_3\\
&&\quad +(q-1)(q-2)(q-3)({\mathcal C}_4+{\mathcal C}_5)\\
{\mathcal C}_3{\mathcal C}_4={\mathcal C}_4{\mathcal C}_3&=&2(q-1)(q-3){\mathcal C}_3+(q-1)(q-2)({\mathcal C}_4+{\mathcal C}_5)\\
{\mathcal C}_3{\mathcal C}_5={\mathcal C}_5{\mathcal C}_3&=&2(q-1)(q-3){\mathcal C}_3+(q-1)(q-2)({\mathcal C}_4+{\mathcal C}_5)\\
{\mathcal C}_4^2&=&2(q-1)^2({\mathcal C}_0+{\mathcal C}_1)+2(q-1){\mathcal C}_3\\
&&\quad +(q-1)(q-2){\mathcal C}_4+(q-1){\mathcal C}_5\\
{\mathcal C}_4{\mathcal C}_5={\mathcal C}_5{\mathcal C}_4&=&2(q-1)^2{\mathcal C}_2+2(q-1){\mathcal C}_3\\
&&\quad +(q-1){\mathcal C}_4+(q-1)(q-2){\mathcal C}_5\\
{\mathcal C}_5^2&=&2(q-1)^2({\mathcal C}_0+{\mathcal C}_1)+2(q-1){\mathcal C}_3\\
&&\quad +(q-1)(q-2){\mathcal C}_4+(q-1){\mathcal C}_5
\end{eqnarray*}
The reader should be warned that ${\mathcal C}_1$ behaves not exactly like $(q-2){\mathcal C}_0$.

We leave it to the reader to write down matrices for multiplications 
with basis-elements and to compute the complete list of minimal idempotents.

Additional association schemes are given by 
${\mathcal C}_0,{\mathcal C}_1,{\mathcal C}_2,{\mathcal C}_3,{\mathcal C}_4+{\mathcal C}_5$ and ${\mathcal C}_0,{\mathcal C}_1+{\mathcal C}_2,{\mathcal C}_3,{\mathcal C}_4+{\mathcal C}_5$.
It is also possible to split ${\mathcal C}_1$ and/or ${\mathcal C}_2$ according to subgroups
of $\mathbb F_q^*$ into several matrices (or classes, as they are sometimes
called). 

We discuss now with a little bit more details the smallest 
interesting association scheme with classes ${\tilde{\mathcal C}}_0,{\tilde{\mathcal C}}_1={\mathcal C}_1+{\mathcal C}_2,{\tilde{\mathcal C}}_2={\mathcal C}_3,{\tilde{\mathcal C}}_3={\mathcal C}_4+{\mathcal C}_5$. Products are given by
${\tilde{\mathcal C}}_0X=X{\tilde{\mathcal C}}_0$ and
\begin{eqnarray*}
{\tilde{\mathcal C}}_1^2&=&(2q-3){\tilde{\mathcal C}}_0+2(q-2){\tilde{\mathcal C}}_1\\
{\tilde{\mathcal C}}_1{\tilde{\mathcal C}}_2={\tilde{\mathcal C}}_2{\tilde{\mathcal C}}_1&=&(2q-3){\tilde{\mathcal C}}_2\\
{\tilde{\mathcal C}}_1{\tilde{\mathcal C}}_3={\tilde{\mathcal C}}_3{\tilde{\mathcal C}}_1&=&(2q-3){\tilde{\mathcal C}}_3\\
{\tilde{\mathcal C}}_2^2&=&(q-1)^2(q-2)({\tilde{\mathcal C}}_0+{\tilde{\mathcal C}}_1)+(q-1)(q-3)(q-4){\tilde{\mathcal C}}_2\\
&&\quad +(q-1)(q-2)(q-3){\tilde{\mathcal C}}_3\\
{\tilde{\mathcal C}}_2{\tilde{\mathcal C}}_3={\tilde{\mathcal C}}_3{\tilde{\mathcal C}}_2&=&4(q-1)(q-3){\tilde{\mathcal C}}_2+2(q-1)(q-2){\tilde{\mathcal C}}_3\\
{\tilde{\mathcal C}}_3^2&=&4(q-1)^2({\tilde{\mathcal C}}_0+{\tilde{\mathcal C}}_1)+8(q-1){\tilde{\mathcal C}}_2+2(q-1)^2{\tilde{\mathcal C}}_3
\end{eqnarray*}
Matrices $M_0,\dots,M_3$ corresponding to multiplication by ${\tilde{\mathcal C}}_0,\dots,{\tilde{\mathcal C}}_3$
are given by
$$
{\tilde{\mathcal C}}_0=\left(\begin{array}{cccc}1&0&0&0\\0&1&0&0\\0&0&1&0\\0&0&0&1\end{array}\right),\quad {\tilde{\mathcal C}}_1=\left(\begin{array}{cccc}0&2q-3&0&0\\1&2(q-2)&0&0\\0&0&2q-3&0\\0&0&0&2q-3\end{array}\right)$$
$${\tilde{\mathcal C}}_2=\left(\begin{array}{cccc}0&0&(q-1)^2(q-2)&0\\
0&0&(q-1)^2(q-2)&0\\1&2q-3&(q-1)(q-3)(q-4)&4(q-1)(q-3)\\
0&0&(q-1)(q-2)(q-3)&2(q-1)(q-2)\end{array}\right),$$
$${\tilde{\mathcal C}}_3=\left(\begin{array}{cccc}0&0&0&4(q-1)^2\\
0&0&0&4(q-1)^2\\0&0&4(q-1)(q-3)&8(q-1)\\
1&2q-3&2(q-1)(q-2)&2(q-1)^2\end{array}\right)$$ 
and minimal idempotents are given by
\begin{eqnarray*}
\beta_0&=&\frac{1}{q^3-q}({\tilde{\mathcal C}}_0+{\tilde{\mathcal C}}_1+{\tilde{\mathcal C}}_2+{\tilde{\mathcal C}}_3),\\
\beta_1&=&\frac{q-2}{2(q^2-q)}({\tilde{\mathcal C}}_0+{\tilde{\mathcal C}}_1)+\frac{1}{q(q-1)^2}{\tilde{\mathcal C}}_2-\frac{q-2}{2q(q-1)^2}{\tilde{\mathcal C}}_3,\\
\beta_2&=&\frac{2q-3}{2(q-1)}{\tilde{\mathcal C}}_0-\frac{1}{2(q-1)}{\tilde{\mathcal C}}_1,\\
\beta_3&=&\frac{1}{q^2-1}({\tilde{\mathcal C}}_0+{\tilde{\mathcal C}}_1)-\frac{2}{(q-1)^2(q+1)}{\tilde{\mathcal C}}_2+\frac{q-3}{2(q-1)^2(q+1)}{\tilde{\mathcal C}}_3.
\end{eqnarray*}
The coefficient of ${\tilde{\mathcal C}}_0$ multiplied by $(q^3-q)$ gives the dimension of the associated eigenspace. Eigenvalues of (with multiplicities) of generators, obtained
by evaluating the characters $\beta_0,\dots,\beta_4$ on ${\tilde{\mathcal C}}_0,\dots,{\tilde{\mathcal C}}_3$,
are given by
$$\begin{array}{c||c|cccc}
&\hbox{dim}&{\tilde{\mathcal C}}_0&{\tilde{\mathcal C}}_1&{\tilde{\mathcal C}}_2&{\tilde{\mathcal C}}_3\\
\hline
\beta_0&1&1&2q-3&(q-1)^2(q-2)&4(q-1)^2\\
\beta_1&\frac{(q+1)(q-2)}{2}&1&2q-3&2(q-1)&4(1-q)\\
\beta_2&\frac{q(2q-3)(q+1)}{2}&1&-1&0&0\\
\beta_3&q&1&2q-3&2(1-q)(q-2)&2(q-1)(q-3)
\end{array}$$

\begin{rem} There exists a few more exotic variation 
of this construction. An example is given by
partitioning the elements of
$\mathrm{SL}_2(\mathbb F_5)$ according to 
the $30$ possible values of the Legendre symbol $\left(x\over 5\right)$
on entries.
Since $-1$ is a square modulo $5$,
these classes are well-defined on $\mathrm{PSL}_2(\mathbb F_5)$
which is isomorphic to the simple group $A_5$.
Details (and a few similar examples) 
will hopefully appear in a future paper.
\end{rem}


\section{Representation-theoretic aspects}\label{sectrepr}

In this Section, we work over $\mathbb C$ for the sake of simplicity.

\subsection{Traces}
Left-multiplication
by the identity $\frac{1}{q-1}A$ of $\mathcal{SC}_{\mathbb C}$ defines an idempotent 
on $\mathbb Q[\mathrm{SL}_2(\mathbb F_q)]$
whose trace is the dimension $\frac{1}{q-1}(q^3-q)=q(q+1)$ of the 
non-trivial eigenspace. Indeed,
every non-trivial element of $\mathrm{SL}_2(\mathbb F_q)$ has trace 
$0$ and the identity-matrix has trace $q^3-q$ since it
fixes all $q^3-q$ elements of $\mathrm{SL}_2(\mathbb F_q)$.
A basis of the non-trivial $q(q+1)$-dimensional 
eigenspace of $\frac{1}{q-1}A$ is given by sums over all 
matrices with rows representing two 
distinct fixed elements of the projective line over $\mathbb F_q$.

The traces $\mathrm{tr}(\pi_1),\dots,\mathrm{tr}(\pi_4)$ of the four
minimal central projectors $\pi_0,\dots,\pi_4$ of $\mathcal{SC}_{\mathbb C}$
are equal to $q^3-q$ times the coefficient of $A$ in $\pi_i$.
They are thus given by
$$\begin{array}{c|cccc}
&\pi_1&\pi_2&\pi_3&\pi_4\\
\hline
\mathrm{trace}&1&\frac{(q+1)(q-2)}{2}&\frac{q(q-1)}{2}&2q
\end{array}$$
and we have 
$$q(q+1)=\mathrm{tr}(\pi_1)+\mathrm{tr}(\pi_2)+\mathrm{tr}(\pi_3)+
\mathrm{tr}(\pi_4),$$
as expected.

\subsection{Characters}

Since simple matrix-algebras of $\mathbb C[\mathrm{SL}_2(\mathbb F_q)]$
are indexed by characters of $\mathbb C[\mathrm{SL}_2(\mathbb F_q)]$, 
it is perhaps interesting to understand all irreducible 
characters involved in idempotents of $\mathcal{SC}_\mathbb C\subset
\mathbb C[\mathrm{SL}_2(\mathbb F_q)]$.

The algebra $\mathcal{SC}_{\mathbb C}$ is in some sense almost ``orthogonal''
to the center of $\mathbb C[\mathrm{SL}_2(\mathbb F_q)]$. The algebra
$\mathcal{SC}_{\mathbb C}$ should thus involve many
different irreducible characters of $\mathrm{SL}_2(\mathbb F_q)$.
We will see that this is indeed the case.

We decompose first the identity 
$\frac{1}{q-1}A$ according 
to irreducible characters of $\mathbb C[\mathrm{SL}_2(\mathbb F_q)]$. 
We refine this 
decomposition to the minimal central idempotents $\pi_1,\dots,\pi_4$ 
of $\mathcal{SC}_\mathbb C$ in Section \ref{subsectdecomppiiqodd}.

For simplicity we work over $\mathrm{GL}_2(\mathbb F_q)$ which has 
essentially the same character-theory as $\mathrm{SL}_2(\mathbb F_q)$.
We work over $\mathbb C$ and we identify (irreducible) characters
with the corresponding (irreducible) representations.

In order to do this, we introduce $F=\sum_{\lambda\in\mathbb F_q^*}\left(\begin{array}{cc}
1&0\\0&\lambda\end{array}\right)\in \mathbb Z[\mathrm{GL}_2(\mathbb F_q)]$.
We have $F^2=(q-1)F$ and $FX=XF$ for $X\in\{A,B,C,D_\pm, E_\pm\}$, considered
as an element of $\mathbb Z[\mathrm{GL}_2(\mathbb F_q)]$.
The map $X\longmapsto \frac{1}{q-1}FX$ preserves traces and 
defines an injective 
homomorphisme of $\mathcal{SC}$ into $\mathbb Q[\mathrm{GL}_2(\mathbb F_q)]$.

We use the conventions of Chapter 5 of \cite{FH} for
conjugacy classes of $\mathrm{GL}_2(\mathbb F_q)$. More precisely, 
we denote by $a_x$ conjugacy classes of central diagonal matrices 
with common diagonal value $x$ in $\mathbb F_q^*$, by $b_x$ 
conjugacy classes given by 
multiplying unipotent matrices by a scalar $x$ in $\mathbb F_q^*$, by $c_{x,y}$
conjugacy classes with two distinct eigenvalues $x,y\in\mathbb F_q^*$
and by $d_\xi$ conjugacy classes with two conjugate eigenvalues $\xi,\xi^q\in
\mathbb F_{q^2}^*\setminus\mathbb F_q$. The number of conjugacy classes of each type is given by
$$\begin{array}{c|c|c|c}
a_x&b_x&c_{x,y}&d_\xi\\
\hline
q-1&q-1&\frac{(q-1)(q-2)}{2}&\frac{q(q-1)}{2}
\end{array}.$$

The character table of 
$\mathrm{GL}_2(\mathbb F_q)$, copied from \cite{FH},
is now given by
$$\begin{array}{r|cccc}
&1&q^2-1&q^2+q&q^2-q\\
&a_x&b_x&c_{x,y}&d_\xi\\
\hline
U_\alpha&\alpha(x^2)&\alpha(x^2)&\alpha(xy)&\alpha(\xi^{q+1})\\
V_\alpha&q\alpha(x^2)&0&\alpha(xy)&-\alpha(\xi^{q+1})\\
W_{\alpha,\beta}&(q+1)\alpha(x)\beta(x)&\alpha(x)\beta(x)&\alpha(x)\beta(y)+
\alpha(y)\beta(x)&0\\
X_\varphi&(q-1)\varphi(x)&-\varphi(x)&0&-(\varphi(\xi)+\varphi(\xi^q))
\end{array}$$
where $\alpha,\beta$ are distinct characters of $\mathbb F_q^*$ and where 
$\varphi$ is a character of $\mathbb F_{q^2}^*$ with $\varphi^{q-1}$ non-trivial. 
The first row indicates the number of elements in a conjugacy class
indicated by the second row. The remaining rows give the character-table.
$U_\alpha$ are the one-dimensional representations factoring through the
determinant. $V_\alpha=V\otimes U_\alpha$ are obtained from the permutation-representation
$V=V_1+U_1$ describing the permutation-action of $\mathrm{GL}_2(\mathbb F_q)$ 
on all $p+1$ points of the projective line over $\mathbb F_q$.
$W_{\alpha,\beta}$ (isomorphic to $W_{\beta,\alpha}$) 
are induced from non-trivial
$1$-dimensional representations of a Borel subgroup (given for example
by all upper triangular matrices). $X_\varphi$ (isomorphic to 
$X_{\overline{\varphi}}$) are the remaining
irreducible representations 
$V_1\otimes W_{\varphi\vert_{\mathbb F_q^*},1}-
W_{\varphi\vert_{\mathbb F_q^*},1}-\mathrm{Ind}_{\varphi}$
with $\mathrm{Ind}_{\varphi}$ obtained by inducing a $1$-dimensional representation $\varphi\not=\varphi^q$
of a cyclic subgroup isomorphic to $\mathbb F_{q^2}^*$.

The trace of the idempotent $\pi=\frac{1}{(q-1)^2}FA$ in irreducible 
representations of $\mathrm{GL}_2(\mathbb F_q)$ is now given by:
\begin{eqnarray*}
U_\alpha&:&\frac{1}{(q-1)^2}\left(\sum_{x\in \mathbb F_q^*}\alpha(x)\right)^2\\
V_\alpha&:&q\left(\frac{1}{q-1}(\sum_{x\in \mathbb F_q^*}\alpha(x^2)+
\frac{1}{(q-1)^2}\left(\sum_{x\in \mathbb F_q^*}\alpha(x)\right)^2\right)\\
W_{\alpha,\beta}&:&(q+1)\left(\frac{1}{q-1}\sum_{x\in \mathbb F_q^*}\alpha(x)\beta(x)+
\frac{2}{(q-1)^2}\left(\sum_{x\in \mathbb F_q^*}\alpha(x)\right)\left(\sum_{x\in \mathbb F_q^*}\beta(x)\right)\right)\\
X_\varphi&:&(q-1)\sum_{x\in\mathbb F_q^*}\varphi(x)
\end{eqnarray*}
The factors $q,q\pm 1$ are the multiplicities of the irreducible 
representations $V_\alpha,W_{\alpha,\beta}$ and $X_\varphi$ in the
regular (left or right) representation. They are of course equal
to the dimensions of $V_\alpha,W_{\alpha,\beta}$ and $X_\varphi$.

Irreducible representations $U_\alpha$ are involved in $\pi$ only
if $\alpha$ is the trivial character. This corresponds of course 
to the central idempotent $\pi_1$ of $\mathcal{SC}_{\mathbb Q}$.

For $V_\alpha$ we get $2q$ for $\alpha$ trivial, $q$ for odd $q$ 
if $\alpha$ is the quadratic character (defined by the Legendre-symbol
and existing only for odd $q$) 
and $0$ otherwise.

For $W_{\alpha,\beta}$ we get $q+1$ if $\beta=\overline{\alpha},\beta\not=\alpha$
and $0$ otherwise. There are $\frac{q-3}{2}$ such representations for odd $q$ and $\frac{q-2}{2}$ such representations for even $q$.

For $X_\varphi$ we get $q-1$ if the character $\varphi$ of $\mathbb F_{q^2}^*$
with non-trivial $\varphi^{q-1}$ 
restricts to the trivial character of 
$\mathbb F_q^*$  and $0$ otherwise. Characters of $\mathbb F_{q^2}^*$
with trivial restrictions to $\mathbb F_q^*$ are in one-to-one correspondence
with characters of the additive group $\mathbb Z/(q+1)\mathbb Z$.
The character $\varphi^{q-1}$ is trivial for two 
of them (the trivial and the quadratic one) if $q$ is odd.
This gives thus $\frac{q-1}{2}$ such irreducible representations for odd $q$.
For even $q$ we get $\frac{q}{2}$ irreducible representations.

All irreducible representations of $\mathrm{GL}_2(\mathbb F_q)$
involved in $\pi=\frac{1}{(q-1)^2}FA$ have irreducible restrictions
to $\mathrm{SL}_2(\mathbb F_q)$. 

The following table sums up contributions of all different irreducible characters to the trace of $\pi=\frac{1}{(q-1)^2}FA$:
$$\begin{array}{c||c|c|c|c}
q&U&V&W&X\\
\hline
\hbox{odd}&1&3q&(q+1)\frac{q-3}{2}&(q-1)\frac{q-1}{2}\\
\hbox{even}&1&2q&(q+1)\frac{q-2}{2}&(q-1)\frac{q}{2}\\
\end{array}$$
and we have 
$$\mathrm{tr}(\pi)=q(q+1)=1+3q+(q+1)\frac{q-3}{2}+(q-1)\frac{q-1}{2},$$
respectively
$$\mathrm{tr}(\pi)=q(q+1)=1+2q+(q+1)\frac{q-2}{2}+(q-1)\frac{q}{2},$$
as expected.

\subsection{Decompositions of $\pi_1,\dots,\pi_4$ for even $q$}\label{subsectdecomppiiqodd}

Over a finite field of characteristic $2$, conjugacy classes are involved in
generators of $\mathcal {SC}$ as follows:
$$\begin{array}{r|ccccccc}
&q-1&q-1&\frac{(q-1)(q-2)}{2}&\frac{q(q-1)}{2}\\
&a_x&b_x&c_{x,y\not=x}&d_{\xi}\\
\hline
FA&1&0&2&0\\
FB&0&q-1&0&0\\
FC&0&(q-1)(q-2)&(q-1)(q-4)&(q-1)(q-2)\\
FD_\pm&0&q-1&2(q-1)&0\\
FE\pm&0&0&q-1&(q-1)\\
\hline
&1&q^2-1&q^2+q&q^2-q\\
\end{array}$$

This implies the following decompositions of central idempotents of
$\mathcal{SC}$: The idempotent $\pi_1$ of rank $1$ is involved with multiplicity
$1$ in the trivial representation of type $U$. 
The idempotent $\pi_2$ of rank $\frac{(q+1)(q-2)}{2}$ (in $\mathbb C[\mathrm{SL}_2(\mathbb F_q)]$) is involved with 
multiplicity $q+1$ in all $\frac{q-2}{2}$ relevant characters of type W.
The idempotent $\pi_3$ of rank $\frac{q(q-1)}{2}$ is involved with multiplicity
$q-1$ in all $\frac{q}{2}$ relevant characters of type . Finally, the idempotent $\pi_4$
is involved with multiplicity $q$ in the irreducible representation $V$.

\subsection{Decompositions of $\pi_1,\dots,\pi_4$ for odd $q$}\label{subsectdecomppiiqodd}

Over a finite field of odd characteristic, contributions of the different conjugacy 
classes to the elements $FX$ (for $X$ a generator of $\mathcal{SC}$) are given by 
$${\scriptsize \begin{array}{r|ccccccc}
&q-1&q-1&\frac{q-1}{2}&\frac{(q-1)(q-3)}{2}&\frac{q-1}{2}&\frac{(q-1)^2}{2}\\
&a_x&b_x&c_{x,-x}&c_{x,y\not=\pm x}&d_{\xi,\mathrm{tr}(\xi)=0}&d_{\xi,\mathrm{tr}(\xi)\not=0}\\
\hline
FA&1&0&2&2&0&0\\
FB&0&0&q-1&0&q-1&0\\
FC&0&(q-1)(q-3)&(q-1)(q-3)&(q-1)(q-4)&(q-1)^2&(q-1)(q-2)\\
FD_\pm&0&q-1&2(q-1)&2(q-1)&0&0\\
FE\pm&0&q-1&0&q-1&0&q-1\\
\hline
&1&q^2-1&q^2+q&q^2+q&q^2-q\\
\end{array}}$$
and conjugacy classes of $\mathrm{GL}_2(\mathbb F_q)$ are involved in 
projectors $\tilde\pi_1=\frac{1}{q-1}\pi_1F,\dots,\tilde\pi_4=\frac{1}{q-1}\pi_4F$ with the following coefficients
$$\begin{array}{r|cccccc|c}
&q-1&q-1&\frac{q-1}{2}&\frac{(q-1)(q-3)}{2}&\frac{q-1}{2}&\frac{(q-1)^2}{2}&\\
&a_x&b_x&c_{x,-x}&c_{x,y\not=\pm x}&d_{\xi,\mathrm{tr}(\xi)=0}&d_{\xi,\mathrm{tr}(\xi)\not=0}&\\
\hline
\tilde\pi_1&\frac{1}{(q-1)^2q(q+1)}&\frac{1}{q(q-1)}&\frac{1}{(q-1)^2}&\frac{1}{(q-1)^2}&\frac{1}{q^2-1}&\frac{1}{q^2-1}\\
\tilde\pi_2&\frac{q-2}{2q(q-1)^2}&\frac{-1}{q(q-1)}&\frac{q-3}{2(q-1)^2}&\frac{-1}{(q-1)^2}&\frac{1}{2(q-1)}&0\\
\tilde\pi_3&\frac{1}{2(q^2-1)}&0&\frac{-1}{2(q-1)}&0&\frac{-1}{2(q+1)}&\frac{1}{q^2-1}\\
\tilde\pi_4&\frac{2}{(q-1)^2(q+1)}&0&\frac{2}{(q-1)^2}&\frac{2}{(q-1)^2}&\frac{-2}{q^2-1}&\frac{-2}{q^2-1}\\
\end{array}$$

The idempotent $\tilde\pi_1$ (or equivalently the idempotent $\pi_1$ of 
$\mathcal{SC}$) is only involved in the trivial representation with multiplicity $1$. The idempotent $\tilde\pi_4$ appears (with multiplicity $q$) only in the
unique non-trivial $q$-dimensional irreducible representation $V$ involved 
in the permutation-representation of $\mathrm{GL}_2(\mathbb F_q)$ acting on 
all $q+1$ points of the projective line over $\mathbb F_q$.

The character $V_L=V\otimes \alpha_L$ where $\alpha_L$ is the quadratic 
character of $\mathbb F_q$ given by the Legendre-symbol has mean-values given
by 
$$\begin{array}{r|cccccc|c}
&q-1&q-1&\frac{q-1}{2}&\frac{(q-1)(q-3)}{2}&\frac{q-1}{2}&\frac{(q-1)^2}{2}&\\
&a_x&b_x&c_{x,-x}&c_{x,y\not=\pm x}&d_{\xi,\mathrm{tr}(\xi)=0}&d_{\xi,\mathrm{tr}(\xi)\not=0}&\\
\hline
V_L&q&0&\left(\frac{-1}{q}\right)&-\frac{1+\left(\frac{-1}{q}\right)}{q-3}&\left(\frac{-1}{q}\right)&\frac{1-\left(\frac{-1}{q}\right)}{q-1}
\end{array}$$
on conjugacy-classes. The character $V_L$ involves thus $\tilde\pi_2$
(or $\pi_2$) with multiplicity $q$ if $q\equiv 1\pmod 4$ and 
$\tilde\pi_3$
(or $\pi_3$) with multiplicity $q$ if $q\equiv 3\pmod 4$.

Mean-values of a character $W_{\alpha,\beta}$ with non-real 
$\beta=\overline{\alpha}$ on conjugacy classes are given by
$$\begin{array}{r|cccccc|c}
&q-1&q-1&\frac{q-1}{2}&\frac{(q-1)(q-3)}{2}&\frac{q-1}{2}&\frac{(q-1)^2}{2}&\\
&a_x&b_x&c_{x,-x}&c_{x,y\not=\pm x}&d_{\xi,\mathrm{tr}(\xi)=0}&d_{\xi,\mathrm{tr}(\xi)\not=0}&\\
\hline
W_{\alpha,\beta}&q+1&1&2\alpha(-1)&\frac{-2(1+\alpha(-1))}{q-3}&0&0
\end{array}$$
They depend on the value $\alpha(-1)\in\{\pm 1\}$. For $q\equiv 1\pmod 4$ there are $\frac{q-5}{4}$ such characters $W_{\alpha,\beta}$ with $\alpha(-1)=1$
and $\frac{q-1}{4}$ such characters with $\alpha(-1)=-1$.
For $q\equiv 3\pmod 4$, there are the same number $\frac{q-3}{4}$
of such characters for both possible values of $\alpha(-1)$.
Such representations involve $\tilde \pi_2$ (or $\pi_2$ of $\mathcal{SC}$)
with multiplicity $q+1$
if $\alpha(-1)=1$ and $\tilde \pi_3$ with multiplicity $q+1$ otherwise.

We consider now a non-real character $\varphi$ of $\mathbb F_{q^2}^*$ with
trivial restriction to $\mathbb F_q^*$. Since $\xi^2$ belongs to
$\mathbb F_q^*$ if $\xi$ has trace $0$, we have 
$\sigma=\varphi(\xi)\in \{\pm 1\}$.

Mean-values on conjugacy-classes 
of a character $X_\varphi$ with $\varphi$ as above are thus given by
$$\begin{array}{r|cccccc|c}
&q-1&q-1&\frac{q-1}{2}&\frac{(q-1)(q-3)}{2}&\frac{q-1}{2}&\frac{(q-1)^2}{2}&\\
&a_x&b_x&c_{x,-x}&c_{x,y\not=\pm x}&d_{\xi,\mathrm{tr}(\xi)=0}&d_{\xi,\mathrm{tr}(\xi)\not=0}&\\
\hline
X_{\varphi}&q-1&-1&0&0&-2\sigma&\frac{2(1+\sigma)}{q-1}
\end{array}$$
For $q\equiv 1\pmod 4$, there are $\frac{q-1}{4}$ such characters for both possible values of $\sigma$. For $q\equiv 3\pmod 4$, the value $\sigma=1$ is achieved by $\frac{q-3}{4}$ such representations and the value $\sigma=-1$ by 
$\frac{q+1}{4}$ representations.

Each such representation with $\sigma=1$ involves $\tilde\pi_3$ (or $\pi_3$)
with multiplicity $q-1$ and each such representation with $\sigma=-1$ involves $\tilde\pi_2$ (or $\pi_2$)
with multiplicity $q-1$.

I thank M. Brion and P. de la Harpe for useful discussions and remarks.

\noindent Roland BACHER, Univ. Grenoble Alpes, Institut 
Fourier (CNRS UMR 5582), 38000 Grenoble, France.

\noindent e-mail: Roland.Bacher@ujf-grenoble.fr


\begin{thebibliography}{99}








\bibitem{B} R.A. Bailey, {\it Association Schemes: Designed Experiments, Algebra and Combinatorics}, Cambridge University Press, 2004.

\bibitem{Ba} E. Bannai, 
Association Schemes and Fusion Algebras, Journal of algebraic Combinatorics 
{\bf 2} (1993), 327--344.

\bibitem{FH} W. Fulton, J. Harris, {\it Representation Theory, A First Course},
Springer, 1991.



\end{thebibliography}
\end{document}